\documentclass[12pt]{iopart}

\expandafter\let\csname equation*\endcsname\relax
\expandafter\let\csname endequation*\endcsname\relax
\usepackage{mathtools}

\usepackage{amsmath}
\usepackage{algorithm,algorithmic}
\usepackage{mathrsfs} 
\usepackage{latexsym}
\usepackage{amsfonts}
\usepackage{changebar}
\usepackage{amsthm}
\usepackage{amssymb}
\usepackage{graphicx}
\usepackage{fullpage}
\usepackage{xcolor}
\usepackage{color}
\usepackage{caption}
\usepackage{subcaption}
\usepackage{listings}
\usepackage{hyperref}
\usepackage{tikz}
\usepackage{graphicx,scalerel}
\usepackage{placeins}

\newcommand{\R}{\mathbb{R}}
\renewcommand{\S}{\mathbb{S}}
\newcommand{\argmin}[1]{\underset{{#1}}{\operatorname{argmin}}}
\renewcommand{\mat}[1]{\mathbf{{#1}}}
\renewcommand{\vec}[1] {\ensuremath{\boldsymbol{#1}}}
\newcommand{\Gamman}{\mat{\Gamma}_{\text{noise}}}

\newcommand{\mpr}{m_{\text{pr}}}
\newcommand{\mupost}{\mu_{\text{post}}}
\newcommand{\mupr}{\mu_{\text{pr}}}

\newcommand{\m}{\vec{m}}

\newcommand{\thth}{\vec{\theta}}
\newcommand{\Cpr}{\mathcal{C}_{\text{pr}}}

\newcommand{\likei}{\pi_{\text{like}}(\vec{y}|m;\thth)}
\newcommand{\like}{\pi_{\text{like}}(\vec{y}|\vec{m};\thth)}

\newcommand{\nDR}{D^{\text{R}}}
\newcommand{\DR}{\mat{D}^{\text{R}}}

\newcommand{\DM}{\mat{D}^{\text{M}}}
\newcommand{\DMi}{\mat{D}^{\text{M}_i}}

\newcommand{\mip}[2]{\langle {#1},{#2}\rangle_{\mat{M}}}

\begin{document}
\title{Hyper-differential sensitivity analysis for nonlinear Bayesian inverse 
problems} 
\author{Isaac Sunseri \dag, Alen Alexanderian \dag ,
 Joseph Hart \ddag, \\ Bart van Bloemen Waanders \ddag }
\address{\dag\ North Carolina State University}
\address{\ddag\ Sandia National Laboratories}
\date{\today}

\begin{abstract}
We consider hyper-differential sensitivity analysis (HDSA) of nonlinear Bayesian inverse
problems governed by PDEs with infinite-dimensional parameters. In previous
works, HDSA has been used to assess the sensitivity of the solution of 
deterministic inverse problems to additional model uncertainties and also 
different types of measurement data. In the present work, we extend HDSA 
to the class of Bayesian inverse problems governed by PDEs. The focus is
on assessing the sensitivity of certain key quantities derived from the 
posterior distribution. Specifically, we focus on analyzing the sensitivity of  
the MAP point and the Bayes risk and make full use of the information
embedded in the Bayesian inverse problem.
After establishing our mathematical framework for HDSA of Bayesian inverse 
problems, we present a detailed computational approach for computing the
proposed HDSA indices. We examine the effectiveness of the proposed approach 
on a model inverse problem governed by a PDE for heat conduction.
\end{abstract}


\section{Introduction} \label{sec:Introduction} 
Many natural phenomena can be described by systems of partial differential
equations (PDEs).  The governing PDEs, however, often include parameters that are unknown and
challenging to measure directly. This gives rise to inverse problems, in which one uses
the PDE model and measurement data to estimate the unknown model
parameters. 
In this article we consider 
Bayesian inverse problems~\cite{Tarantola05,Stuart10}, whose solution
is 
a posterior distribution
that is informed by both our prior knowledge and the data measurements. 
Specifically, we focus on Bayesian 
inverse problems governed by PDEs with infinite-dimensional
parameters.

In addition to the parameters being estimated, the governing PDEs typically contain parameters that are uncertain but needed for a full
model specification. 
For clarity, we refer to the parameters being estimated
by the inverse problem as \textit{inversion parameters} and call the additional
model parameters the \emph{auxiliary parameters}.  Another source of
uncertainty in the inverse problem arises from the parameters specifying the
experimental conditions, such as the location of measurement devices or their
accuracy. We call these the \textit{experimental parameters}.  Throughout the
article we will refer to the union of auxiliary and experimental parameters as
\textit{complementary parameters}. Our goal in this article is to develop
methods for assessing the sensitivity of the solution of a Bayesian inverse
problems with respect to perturbations of complementary parameters; 
see Section~\ref{sec:MotivatingEx} for a simple illustrative example.

Understanding the sensitivity of an inverse problem to
complementary parameters is important.  These parameters may differ from their
measured or estimated values, which in turn will result in a solution different
from the one we would obtain if we had access to perfect measurements and true
auxiliary parameters. Determining the sensitivity of the solution to
perturbations in these parameters can inform our modeling assumptions and
experimental design practices.  Specifically, this can guide a goal oriented
prioritization of resources and focus efforts on obtaining accurate values for
the important auxiliary parameters. Moreover, if one has data that is
informative to the important auxiliary parameters, the inverse problem may be
redesigned to include these parameters in the set of inversion
parameters.

The present work builds on previous efforts in hyper-differential sensitivity
analysis (HDSA) \cite{SunseriHartEtAl20,
HartvanBloemenWaanders20,Brandes06,Griesse_AD,Griesse_SISC,Griesse_Thesis,
Griesse_part_1,Griesse_part_2,griesse_constraints,griesse_3d}. Traditional HDSA
uses the derivative of the solution of an optimization problem with respect to
complementary parameters to define sensitivity indices. These indices measure
how much the solution of the optimization problem changes when the
complementary parameters are perturbed. Specifically, in our previous
work~\cite{SunseriHartEtAl20}, which targets deterministic inverse problems, we
define two types of sensitivity indices: pointwise sensitivities, and
generalized sensitivities.  Pointwise indices measure the sensitivity of the  
solution to perturbation in specific complementary
parameters. Generalized indices measure the maximum possible change in the
solution with respect to any unit perturbation of groups of
complementary parameters. These indices 
provide a framework to study the sensitivity
of the inverse problem solution; see~\cite{SunseriHartEtAl20} for more details.

In the present work, we extend HDSA to the class of Bayesian inverse problems,
which allows us to study the change in the posterior distribution with respect
to perturbations of the complementary parameters. We focus on HDSA of Bayesian
inverse problems governed by PDEs with infinite-dimensional parameters.  
Section~\ref{sec:Preliminaries} provides a brief overview of the inverse problems
under study.
HDSA of a Bayesian inverse problem is difficult, because the solution of such
problems is a statistical distribution. In general,
the posterior distribution is difficult to approximate; this 
makes assessing the sensitivity of the posterior to complementary
parameters challenging. A tractable approach is to instead focus on
certain key aspects of the posterior distribution.  
Namely, we consider specific quantities
of interest (QoIs) derived from the posterior distribution to perform HDSA on; we call
such quantities the \emph{HDSA QoIs}.

A first possibility, which we consider in this article, is to assess the
sensitivity of the maximum a posteriori probability (MAP) point to the
complementary parameters. This builds directly on the developments
in~\cite{SunseriHartEtAl20}.  Of greater difficulty is obtaining sensitivities
of a measure of the posterior uncertainty.  A natural setting for defining such
measures is provided by the theory of optimal experimental design
(OED)~\cite{Atkinson92,Pazman86,Pukelsheim93,
ChalonerVerdinelli95,Ucinski05,Alexanderian21}.  Recall that in OED, one seeks
experiments that minimize posterior uncertainty or, more generally,
optimize the statistical quality of the estimated parameters. 
This is done 
by optimizing certain design criteria.  Examples include the A-optimality
criterion, which quantifies the average posterior variance, or the Bayesian
D-optimality criterion, measuring the expected information gain; see
e.g.,~\cite{ChalonerVerdinelli95}. In the present work, we consider the Bayes
risk, which has been used previously in OED for PDE-constrained inverse
problems~\cite{HaberHoreshTenorio_08,HaberHoreshTenorio_10,
HoreshHaberTenorio_10}.  Our motivations for using the Bayes risk as an HDSA
QoI are two-fold. First, the Bayes risk is defined as an average error of
the MAP estimator; see Section~\ref{sec:HDSA_QoIs}.  Thus, 
HDSA of the Bayes risk builds further on
methods for HDSA of the MAP point. Moreover, it is 
well-known~\cite{ChalonerVerdinelli95,AlexanderianGloorGhattas16} that the Bayes
risk, with respect to the $L^2$ loss function, reduces to the A-optimality
criterion, in the case of Gaussian linear Bayesian inverse problem.  Hence, up
to a linearization, the Bayes risk may be considered as a proxy for the average
posterior variance.  Bayes risk is also a common utility function in 
decision theory.  

The contributions of this article are as follows:
\begin{itemize}
\item We develop a mathematical framework to 
assess the sensitivity of the Bayes risk and the MAP point in nonlinear Bayesian inverse problems 
with respect to complementary parameters; see Section~\ref{sec:HDSA}.

\item We present a scalable computational framework
for computing the HDSA indices for the MAP point and the Bayes risk; see
Section~\ref{sec:Computations}.  In that section, we also detail the
computational cost of the various components of the proposed approach. 

\item We present comprehensive numerical results for a model problem of heat
flow across a conductive surface that examine the effectiveness and efficiency
of the proposed approach; see Section~\ref{sec:ModelProblem} for the description of
the model under study and Section~\ref{sec:Results} for our computational
results.

\end{itemize}

\section{An Illustrative Example} \label{sec:MotivatingEx}
\begin{figure}[ht!]
	\centering
	\includegraphics[trim={0 0 0 0},clip,width=0.45\linewidth]{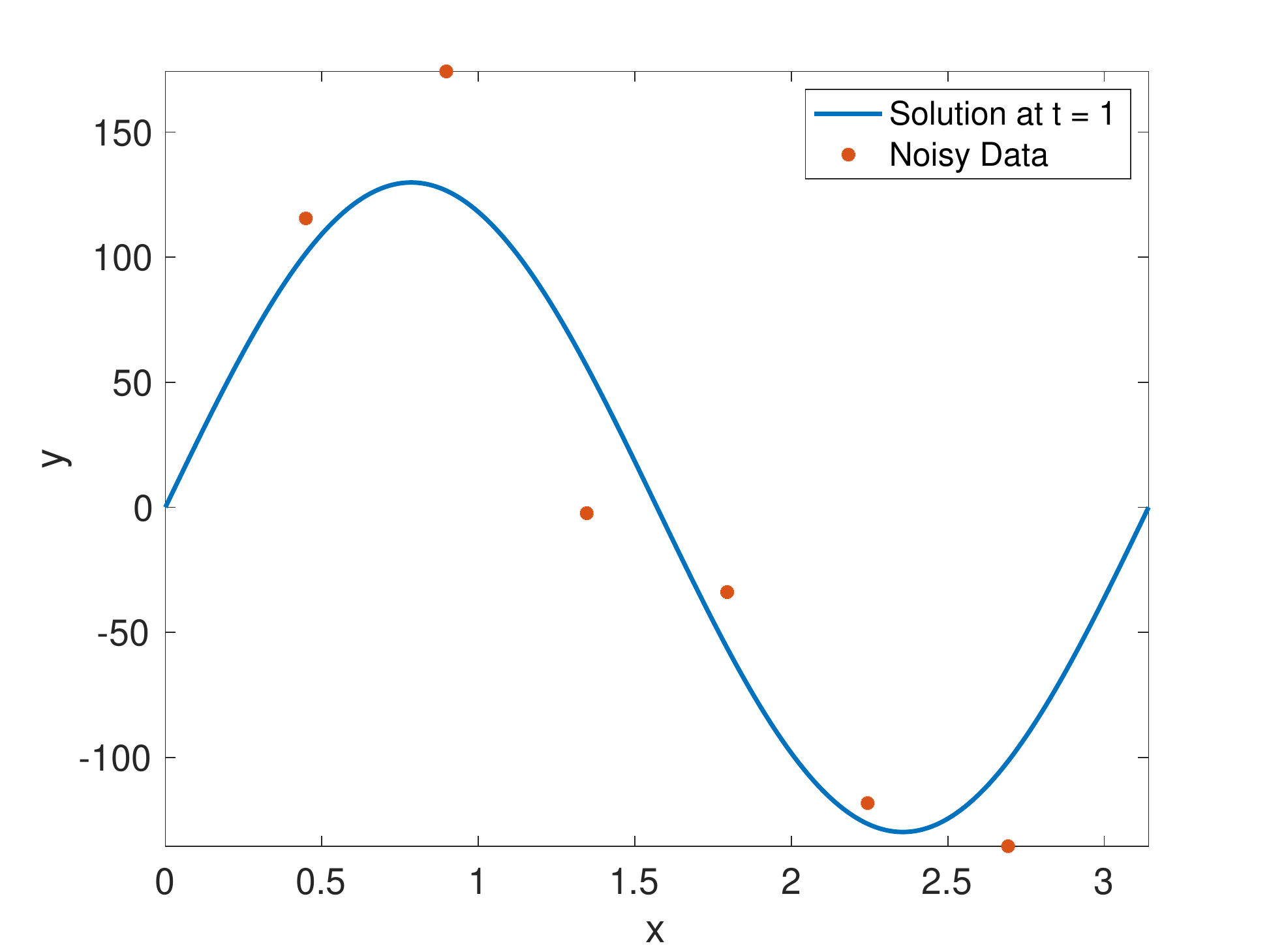}
	\includegraphics[trim={0 0 0 0},clip,width=0.45\linewidth]{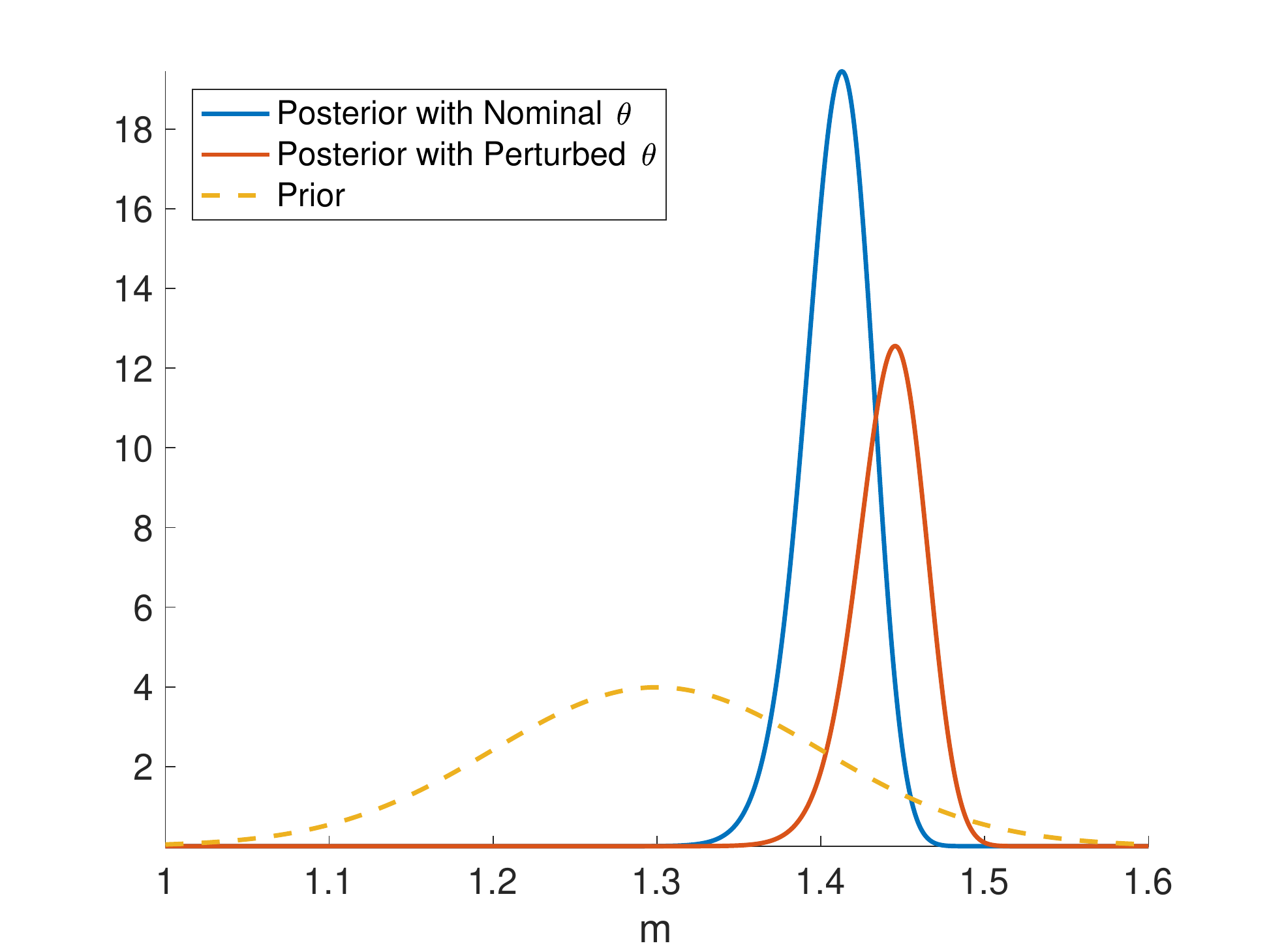}
	\caption{Left: the solution of the PDE at the final time with noisy data measurements, Right: the prior pdf of $m$ and its posterior pdf for both a nominal and perturbed value of the auxiliary parameter $\theta$.}
	\label{fig:examplepdf}
\end{figure}

We consider a simple example to motivate the problem considered in this work, which is to conduct sensitivity analysis on the solution of Bayesian inverse problems. Consider the heat equation, 
\begin{subequations} \label{eq:HeatEx}
	\begin{align} 
	\frac{\partial y(x,t)}{\partial t} &= \exp(m)\frac{\partial^2y(x,t)}{\partial x^2}, &&x \in (0,\pi), t \in (0,1) \\
	y(0,t) &= y(\pi,t) = 0, &&t \in (0,1) \\
	y(x,0) &= \sin(x) + \exp(\theta)\sin(2x), &&x \in (0,\pi). 
	\end{align}
\end{subequations}
In this problem, $m$ is the inversion parameter and $\theta$ is an uncertain auxiliary parameter. For simplicity we let $m$ and $\theta$ be scalars here. Following a Bayesian framework, we endow $m$ with a Gaussian prior $m \sim \mathcal{N}(1.3,.1)$. The problem \eqref{eq:HeatEx} can be solved analytically and the solution is given by,
\[
y(x,t) = \exp(-\exp(m)t)\sin(x) + \exp(-4\theta\exp(m)t)\sin(2x).
\]
We have collected data measurements with additive Gaussian noise at the final time $t = 1$ at six evenly spaced points between $x = 0$ and $x = \pi$ depicted in Figure \ref{fig:examplepdf}. The Gaussian noise is unbiased with a standard deviation of 26 in this example. 
By Bayes rule, the posterior probability density function (pdf) is proportional to the product of the likelihood and prior pdfs. Because the auxiliary parameter is uncertain, we want to understand how the posterior distribution changes as the auxiliary parameter is perturbed. 
We view the posterior pdf of $m$ solved at the nominal value of $\theta = -.3$ and a perturbed value of $\theta = -.29$ in Figure \ref{fig:examplepdf}. 

\FloatBarrier

We notice two primary changes in the posterior pdf as the auxiliary parameter is perturbed. First, there is a shift in the location of the distribution's peak, which equates to a change in the MAP point. Second, there is a change in the spread or variance of the distribution which equates to a changing of the posterior uncertainty. We can see that perturbations of auxiliary parameters can have significant impact upon both of these posterior quantities. 

\section{Bayesian inverse problems and complementary parameters} \label{sec:Preliminaries}

We let $\thth_a$ and $\thth_e$ denote the auxiliary and experimental parameters, and call the augmented parameter
vector $\thth = \begin{bmatrix}
\thth_a \\ \thth_e
\end{bmatrix}$ the complementary parameters.
Note that it is possible to have some auxiliary parameters that are
functions (see~\cite{SunseriHartEtAl20}) but to keep the presentation simple,
we consider finite-dimensional complementary parameters. The precise definition of the experimental parameters is, in general, application dependent, but in what follows we consider a particular case. Our goal
is to provide a comprehensive framework for analyzing the sensitivity of the
solution of the inverse problems under study to perturbations in $\thth$. 

We assume that the governing PDE (the state equation),
represented abstractly by
\begin{equation}
v(u,m,\thth_a) = 0  \label{eq:cont_state},
\end{equation}
has a unique solution $u$ for a given $m$ and fixed auxiliary parameters
$\thth_a$.  The inversion parameter $m$ belongs to an infinite-dimensional
Hilbert space  $\mathcal{M}$ that is equipped with 
an inner product $\langle \cdot,
\cdot \rangle_{\mathcal{M}}$ and the induced norm $\|\cdot\|_{\mathcal{M}}$.
The
state variable $u$ belongs to an infinite dimensional reflexive Banach space
$U$.
In the present work, $\mathcal{M} = L^2(\Omega)$ 
where $\Omega$ is a suitable physical domain and  $\langle \cdot, \cdot \rangle_{\mathcal{M}}$
is the standard $L^2$ inner product.

To infer $m$, we solve a Bayesian inverse problem that uses observed
measurements along with information known about the governing system of PDEs.
We assume that (noisy)
measurement data is related to $m$ according to the following model:
\begin{equation} \label{eq:data_model}
\vec{y} = \vec{F}(m,\thth_a) + \vec{\eta}(\thth_e),
\end{equation}
where $\vec{y}$ is a vector of $n_y$ experimental measurements, $\vec{F}$ the
parameter-to-observable map that takes in the inversion parameter $m$ and maps it to a 
vector of measurements, and $\vec{\eta}$ a vector that models additive Gaussian 
noise, $\vec{\eta} \sim \mathcal{N}(\vec{0},\Gamman(\thth_e))$. Evaluating 
$\vec{F}(m,\thth_a)$ requires 
solving the state equation \eqref{eq:cont_state} followed by application of an observation 
operator $\mathcal{O}$ which evaluates the state $u$ at the $n_y$ sensor locations. 
In the present work we let $\thth_e$ parameterize the 
noise levels of the sensors. This can correspond to situations where the experimental
error at various sensors can be controlled either by repeated measurements or by choice of
the measurement device, or possibly recalibration of existing devices.
Our model for the experimental parameters is detailed in Section \ref{sec:ModelProblem}.

\textbf{The Bayesian inverse problem setup.}
To solve an inverse problem with the data model \eqref{eq:data_model} we 
define the data likelihood pdf $\likei$, which describes the distribution of data 
measurements $\vec{y}$, given a particular inversion parameter $m$.
Given our assumption of an additive Gaussian noise model, we have 
$\vec{y}|m \sim \mathcal{N}(\vec{F}(m),\Gamman)$, and thus, 
\begin{equation} \label{eq:like}
\likei \propto \exp \left( -\frac12(\vec{F}(m,\thth_a)-\vec{y}(\thth_e))^{\top}\Gamman^{-1}(\thth_e)(\vec{F}(m,\thth_a)-\vec{y}(\thth_e)) \right). 
\end{equation}
Note that in \eqref{eq:like} the auxiliary parameters appear in the parameter-to-observable map $\vec{F}$ while we assume the experimental parameters for our problem appear only in the data measurements themselves and the noise covariance matrix. 

In a Bayesian paradigm, we model our uncertainty regarding the inversion
parameter by modeling $m$ as a random variable.  
Accordingly, we endow the inversion parameter $m$ with a prior distribution 
that reflects our 
knowledge of $m$ a priori. In the present work, we 
let the prior distribution law be a Gaussian $\mupr = \mathcal{N}(\mpr,\Cpr)$, with mean $\mpr$ and covariance operator $\Cpr$. We let $\Cpr = \mathcal{A}^{-2}$ where $\mathcal{A}$ is a Laplace-like 
differential operator; see, e.g.,\cite{Alexanderian21,Stuart10, Bui-ThanhGhattasMartinEtAl13}. 
The Gaussian prior measure is meaningful since $\mathcal{A}^{-2}$ is a trace class operator which guarantees bounded variance and almost surely pointwise well-defined samples. 
The prior measure induces the Cameron-Martin space $\mathscr{E} = \text{range}(\Cpr^{1/2})$, 
which is endowed with the following inner product, 
\[
\langle x,y \rangle_{\mathscr{E}} = \langle \mathcal{A}x,\mathcal{A}y \rangle. 
\]
We assume $\mpr \in \mathscr{E}$. 

The definition of the prior measure and the data likelihood completes the
description of the Bayesian inverse problem. 
The solution of this inverse problem is called the
posterior measure $\mupost^{\vec{y}}$, which describes the probability law of
$m$ conditioned on experimental measurements $\vec{y}$. We
will often denote the posterior measure as $\mupost$ for notational simplicity
when no confusion arises from doing so. The Bayes formula takes the following
form in the infinite-dimensional Hilbert
space setting~\cite{Stuart10}
\[
\frac{d\mupost}{d\mupr} \propto \likei.  
\]
Also, for a fixed $\thth$, the maximum a posteriori probability (MAP) estimator of $m$ is found by solving 
\begin{equation} \label{eq:MAP_opt}
m^*(\thth) \coloneqq \argmin{m \in \mathscr{E}} \ 
J(m,\thth), 
\end{equation}
where
\begin{equation}
\label{eq:MAP_cost}
J(m, \thth) \coloneqq 
\frac12 
\big(
\vec{F}(m,\thth_a)-\vec{y}(\thth_e)
\big)^\top
\Gamman^{-1}(\thth_e)\big(\vec{F}(m,\thth_a)-\vec{y}(\thth_e)\big) 
 + \frac12\langle m - \mpr, m-\mpr \rangle_{\mathscr{E}}.
\end{equation}

\textbf{Discretization.}
In the present work, we follow a continuous Galerkin 
finite element discretization
and let $\vec{m}$ and $\vec{u}$ be the discretizations
of their continuous counterparts $m$ and $u$. 
We let $n_m$ be the the dimension of the discretized parameter.
The discretized space is $\R^{n_m}$
equipped with the inner product
\[
\mip{\vec{a}}{\vec{b}} = \vec{a}^\top \mat{M} \vec{b}, 
\quad \vec{a}, \vec{b} \in \R^{n_m},
\] 
where $\mat{M}$ is the finite element mass matrix, and the norm $\| \cdot
\|_\mat{M}$ induced by this inner product.  Note that 
when working with linear operators on $(\R^{n_m}, \| \cdot\|_{\mat{M}})$ 
or linear transformations between $(\R^n, \|\cdot \|_\mat{M})$
and $(\R^n, \| \cdot \|)$, 
where
$\| \cdot \|$ is the Euclidean inner product, the adjoint operators need
to be defined appropriately; see~\cite{Bui-ThanhGhattasMartinEtAl13} 
for this and further details on 
discretization of different components of 
infinite-dimensional
Bayesian inverse problems. In the remainder of this article, we present
the proposed methods in the discretized setting.

\section{HDSA for Nonlinear Bayesian Inverse Problems} \label{sec:HDSA}
In this section, we outline a framework for HDSA of nonlinear Bayesian inverse problems.

\subsection{The HDSA QoIs}\label{sec:HDSA_QoIs}
As discussed in the introduction, we consider two HDSA
QoIs for a Bayesian inverse problem: (i) the MAP point, 
which is obtained by minimizing~\eqref{eq:MAP_cost} and (ii) the Bayes risk.
The Bayes risk, for a fixed vector $\thth$ of complementary parameters, is defined by 
\begin{align} \label{eq:BayesRiskInf}
\Psi_{\text{risk}}(\thth) = \int_{\mathcal{M}} \int_{\mathbb{R}^d} 
\| \vec{m}^*(\thth) - \vec{m} \|_{\mat{M}}^2 \like d\vec{y} \mupr^{n_m}(d\vec{m}). 
\end{align}
Note that here we expressed the Bayes risk for the discretized version of 
the Bayesian inverse problem, and $\vec{m}^*$ is the discretized MAP point. 
The discretized prior measure, which we 
denote by $\mupr^{n_m}$ should be defined appropriately, as described 
in~\cite{Bui-ThanhGhattasMartinEtAl13}.

In practice, Bayes Risk is approximated via sample averaging.
Namely, we draw $n_s$ samples $\{\vec{m}_1,\dots,\vec{m}_{n_s}\}$ 
from the prior distribution to compute 
data samples $\{\vec{y}_1,\dots,\vec{y}_{n_s}\}$ with the forward data model, 
\[
\vec{y_i} = \vec{F}(\vec{m}_i,\thth_a) + \vec{\eta}_i(\thth_e), \qquad i = 1,\dots,n_s,
\]
where $\vec{\eta}_i$ are sample draws from noise distribution 
$\mathcal{N}(0,\Gamman(\thth_e))$. We can then rewrite the approximate 
Bayes Risk as
\begin{align} \label{eq:BayesRiskDis}
\widehat{\Psi}_{\text{risk}}(\thth) =  \frac{1}{n_s} \sum_{i=1}^{n_s} \| \vec{m}^*(\vec{y}_i,\thth) - \vec{m}_i \|_{\mathcal{M}}^2.  
\end{align}

\subsection{Sensitivity operator of Bayes risk}

To compute the sensitivity operator of the approximate Bayes risk estimator, we first discretize \eqref{eq:BayesRiskDis} and denote its discretization by $\widehat{\mat{\Psi}}_{\text{risk}}(\thth)$. We assess the sensitivity of Bayes risk by differentiating $\widehat{\mat{\Psi}}_{\text{risk}}(\thth)$ with respect to $\theta_j$, the $j^{\text{th}}$ component of $\thth$, and evaluate at a set of nominal complementary parameter values $\thth^*$:
\begin{equation} \nonumber
\nDR_j \coloneqq \frac{\partial}{\partial \theta_j} \widehat{\mat{\Psi}}_{\text{risk}}(\thth^*) = \frac{2}{n_s} \sum_{i=1}^{n_s} \frac{\partial}{\partial \theta_j}(\vec{m}^*(\vec{y}_i,\thth^*))^{\top}\mat{M}\vec{m}^*(\vec{y}_i,\thth^*) - \frac{\partial}{\partial \theta_j}(\vec{m}^*(\vec{y}_i,\thth^*))^{\top}\mat{M}\vec{m}_i.
\end{equation} 
We define the discretized sensitivity operator of the approximate
Bayes risk as 
\begin{equation}\label{eq:DR_def}
\DR = 
\begin{bmatrix} 
\nDR_1 & \nDR_2 & \dots & \nDR_{n_\theta}\end{bmatrix},
\end{equation}
where
$n_\theta$ denotes the dimension of the complementary parameter vector.
Note that $\DR \tilde{\thth}$ 
can be interpreted as the sensitivity of the approximate Bayes risk with
respect to a perturbation of the complementary parameters in the direction
$\tilde{\thth}$. 

To compute the derivative of the approximate Bayes risk, we need
$\frac{\partial \vec{m}^*}{\partial \theta}(\vec{y}_i,\thth^*)$, $i = 1, \ldots, n_s$, which
measures the sensitivity of the MAP points (for each data sample $\vec{y}_i$) to the complementary parameters \footnote{As discussed further in 
the Section~\ref{sec:Computations}, we only need to compute the action 
of 
this sensitivity operator on vectors}. 
For clarity, we denote the 
discretized cost functional by $\vec{J}$. 
As discussed in~\cite{SunseriHartEtAl20}, under mild 
regularity assumptions~\cite{HartvanBloemenWaanders20,Brandes06,Ambrosetti95}
using the implicit function theorem, we obtain 
\begin{align} \label{eq:SensOp}
\DM = \frac{\partial \vec{m}^*}{\partial \thth} = -\left(\frac{\partial^2 \vec{J}}{\partial \vec{m}^2}\right)^{-1}\frac{\partial^2 \vec{J}}{\partial \vec{m} \partial \thth} = -\mat{H}^{-1}\mat{B},
\end{align}
where $\mat{H}$ and $\mat{B}$ are evaluated at the solution $\vec{m}^*$ with fixed nominal parameters $\thth^*$. By averaging these computed sensitivities over
the number of data samples $n_s$, we can simultaneously measure both the average MAP point and Bayes risk sensitivities.

It is important to note the significance of this process. In a deterministic formulation ~\cite{SunseriHartEtAl20}, the sensitivities of the inverse problem solution $\vec{m}^*$ require data measurements to compute. That is, some experimental measurements would be needed before conducting sensitivity analysis. In contrast, the method proposed here does not require experimental measurements and can be computed a priori by using the information encoded in the Bayesian inverse problem to generate likely data realizations. This makes the methodology applicable to a broad range of problems where data is not available at the time of performing HDSA.

\subsection{Sensitivity Indices} \label{sec:sens_ind}
Given a sensitivity operator, we define sensitivity indices to provide a scalar which measures the magnitude of the change in the solution with respect to a particular perturbation of the complementary parameters. We first group related complementary parameters together into $K$ subsets. For example, we group data measurements corresponding to the same state variable together, scalar auxiliary parameters form their own group (of size 1), while all parameters defining the discretization of an uncertain function may form another group. Let $\Theta_k$ be the inner product space containing the $k$th set of parameters for $k = 1,\dots,K$ and $\{\vec{b}_k^1,\vec{b}_k^2,\dots,\vec{b}_k^{n_k}\}$ be a basis for $\Theta_k$ of dimension $n_k$. We then define $\{\vec{e}_k^j\}$ as the basis of $\Theta = \Theta_1 \times \Theta_2 \times \dots \times \Theta_K$ for $k = 1,\dots,K$ and $j = 1,\dots,n_k$ where
$$ \vec{e}_k^j = \begin{pmatrix}
\vec{0}_1 & \dots & \vec{0}_{k-1} & \vec{b}_k^j & \vec{0}_{k+1} & \dots & \vec{0}_K
\end{pmatrix}^{\top}.$$

We can now define pointwise sensitivity indices to measure the sensitivity of both the MAP point and the approximate Bayes risk as,
\begin{equation} 
\label{eq:Risk_sens}
S_k^j = \frac{\| \DM \vec{e}_k^j \|_{\mat{M}}}{\| \vec{e}_k^j \|_{\Theta}} \qquad \text{and} \qquad \S_k^j = \frac{| \DR \vec{e}_k^j |}{\| \vec{e}_k^j \|_{\Theta}}
\end{equation}
respectively. These pointwise sensitivities measure the change in the MAP point or Bayes risk, respectively, to a perturbation of the $k$th parameter in the $j$th direction $\vec{b}_k^j$. 

We would also like to determine the importance of the $K$ parameter subgroups relative to one another. To do so, we define generalized sensitivity indices which provide a single measure of sensitivity for each parameter subgroup. Let $\mat{T}_k: \Theta \to \Theta$ be a selection operator that zeros out components of $\thth$ not in $\Theta_k$. We define the generalized sensitivity of the $k$th subgroup of complementary parameters with respect to the MAP point and approximate Bayes risk respectively as, 
\begin{equation}
\label{eq:MAP_sens}
S_k = \max_{\vec{\theta} \in \Theta}
\frac{\|\DM\mat{T}_k\vec{\theta}\|_{\mat{M}}}{\|\vec{\theta}\|_{\Theta}} \qquad \text{and} \qquad 
\S_k = \max_{\vec{\theta} \in \Theta}
\frac{| \DR\mat{T}_k \vec{\theta} |}{\|\vec{\theta}\|_{\Theta}}.
\end{equation}

The generalized sensitivities measure the maximum change that can be observed in the HDSA QoIs to a norm-1 perturbation of the $k$th parameter subgroup. We can interpret this as a ``worst case scenario" sensitivity because it measures the maximum change in the solution. More importantly, the generalized sensitivities provide a single measure of sensitivity for each parameter subgroup which can be used to compare their relative importance, despite their potentially diverse range of physical characteristics. Note that the parameter groupings should be specified by the user and are problem dependent. In the model problem considered in Section $\ref{sec:ModelProblem}$ we allow scalar auxiliary parameters to each consist of their own subgroup while the experimental parameters, corresponding to noise in the data measurements, are grouped together. It is important to note that if a subgroup consists of a single scalar parameter, its pointwise and generalized sensitivities will be identical. We direct the reader to \cite{SunseriHartEtAl20} for additional details on the construction of these sensitivities. 

To compare the MAP point and Bayes risk sensitivities it is important to note that each sensitivity is endowed with specific units. If we were only concerned with a single HDSA QoI, this would not matter because we would be primarily concerned with the relative differences between sensitivities of that measure. When comparing the sensitivities of the MAP point to Bayes risk however, we must normalize with respect to the QoI to compare the sensitivities to each other in a reasonable fashion. To do so, we divide the sensitivities with respect to the MAP point by the average norm of the computed MAP points, $ \frac{1}{n_s}\sum_{i = 1}^{n_s} \| \vec{m}^*(\vec{y}_i, \thth^*) \|_{\mat{M}}$, and the sensitivities with respect to Bayes risk by the computed value of Bayes risk, $\widehat{\mat{\Psi}}_{\text{risk}}(\thth^*)$. 

\section{Computational Methods} \label{sec:Computations}
In this section, we present  computational methods to implement the framework
proposed in Section~\ref{sec:HDSA}. 

\subsection{Computing the sensitivity indices} \label{sec:Computations1}
With the discretization of $\vec{m}$ and its prior, we can write the sensitivity operator
of the approximate Bayes risk with respect to the complementary parameters \eqref{eq:DR_def} as, 
\begin{equation} \label{eq:dis_DR}
\DR = \frac{2}{n_s} \sum_{i=1}^{n_s} (\DMi)^{\top}\mat{M}\left(\vec{m}^*(\vec{y}_i,\thth^*) - \vec{m}_i \right) = \frac{2}{n_s} \sum_{i=1}^{n_s} -\mat{B}_i^{\top}\mat{H}_i^{-\top}\mat{M}\left(\vec{m}^*(\vec{y}_i,\thth^*) - \vec{m}_i \right).
\end{equation}
We let the subscript $i$ on the operators $\DMi\in \R^{n_m \times n_{\theta}}$, $\mat{B}_i \in \R^{n_m \times n_\theta}$, and $\mat{H}_i \in \R^{n_m \times n_m}$ indicate the dependence on the $i$th data sample.

To compute matrix-free actions of $\mat{H}$ and $\mat{B}$ to vectors, we use a discretized formal Lagrangian approach. We note that this method is utilized to both compute sensitivity indices as well as solve for the MAP point. We begin by defining the discrete Lagrangian as 
\begin{equation} \label{eq:disc_Lagrangian}
\mathcal{L}(\vec{u},\vec{m},\vec{p};\thth) = \vec{J}(\vec{m},\thth) - 
\mip{\vec{p}}{\vec{v}(\vec{u},\vec{m},\thth_a)}
\end{equation} 
where $\vec{v}(\vec{u},\vec{m},\thth_a)$ is the discretized form of the PDE $v$, and $\vec{p}$ is the adjoint variable. Next, we use variational derivatives to compute the action of the discretized gradient of the cost function. 
We let $\mathcal{L}_{\vec{p}}[\hat{\vec{p}}]$ denote the variational derivative
of~\eqref{eq:disc_Lagrangian} with respect to $\vec{p}$, acting on
$\hat{\vec{p}}$, with the input arguments suppressed for brevity. A
similar notation is used for the variational derivatives with respect to
$\vec{u}$ and $\vec{m}$. We can also compute the action of the Hessian by
constructing a meta-Lagrangian,
\begin{equation} \label{eq:meta_Lagrangian}
\mathcal{L}^H(\vec{u},\vec{m},\vec{p},\hat{\vec{u}},\hat{\vec{m}},\hat{\vec{p}};\thth) = 
\mathcal{L}_{\vec{p}}[\hat{\vec{p}}] + 
\mathcal{L}_{\vec{u}}[\hat{\vec{u}}] + 
\mathcal{L}_{\vec{m}}[\hat{\vec{m}}]. 
\end{equation}
By computing variational derivatives of the meta-Lagrangian, we can evaluate
the action of the discretized Hessian to vectors. The basic steps of this
solution process are outlined in Algorithm~\ref{alg:H}, and we direct the
reader to \cite{VillaPetraGhattas21, Gunzburger03} for additional details. 

\begin{algorithm}[ht!!]
	\caption{Compute the gradient g($\vec{m}$) and action of the Hessian $\mat{H}$ in the direction $\hat{\vec{m}}$}
	\begin{algorithmic} [1] \label{alg:H}
		\STATE Solve the state equation $\mathcal{L}_{\vec{p}} = 0 \quad$ for the state variable $\vec{u}$
		\STATE Solve the adjoint equation $\mathcal{L}_{\vec{u}} = 0 \quad$ for the adjoint variable $\vec{p}$
		\STATE Evaluate $g(\vec{m})^{\top} = \mathcal{L}_{\vec{m}}$ 
		\STATE Solve the incremental state equation $\mathcal{L}^H_{\vec{p}} = 0 \quad$  for the incremental state variable $\hat{\vec{u}}$
		\STATE Solve the incremental adjoint equation $\mathcal{L}^H_{\vec{u}} = 0 \quad$ for the incremental adjoint variable $\hat{\vec{p}}$
		\STATE Evaluate the Hessian apply $\mat{H}(\vec{m})[\hat{\vec{m}}] = \mathcal{L}^H_{\vec{m}}$ 
	\end{algorithmic}
\end{algorithm}
\FloatBarrier
Next, we discuss computing the action of the mixed derivative operator
$\mat{B}^{\top}$. We follow a similar approach as one used to compute the 
action of the Hessian using the meta-Lagrangian $\mathcal{L}^H$. 
Namely, we differentiate the mata-Lagrangian with respect to 
$\thth$ to obtain
\[
\mathcal{L}^H_{\thth}
(\vec{u},\vec{m},\vec{p},\hat{\vec{u}},\hat{\vec{m}},\hat{\vec{p}};\thth)
[\tilde{\thth}] = \tilde{\thth}^{\top}\mat{B}^{\top} \hat{\m},
\]
where $\hat{\vec{u}}$ and $\hat{\vec{p}}$ satisfy the incremenal state and 
and adjoint equations, respectively.

Note that we can also compute the action of $\mat{B}$ by reversing the order of
differentiation, deriving through the Lagrangian by $\thth$ and the
meta-Lagrangian by $\vec{m}$, which will result in modified incremental
equations. These adjoint based methods provide a computationally efficient
method to evaluate the sensitivity operators $\DM$ and $\DR$. 

To compute the discretized sensitivity operator $\DR$, we must first generate data samples $\vec{y}_i$ for $i = 1,\dots,n_s$ and then evaluate~\eqref{eq:dis_DR} which requires non-trivial computational cost.
We also compute sensitivities of the MAP point, efficiently reusing PDE solves whenever applicable. This process is summarized in Algorithm~\ref{alg:SensInd}. Note that for clarity, we have separated the processes of data generation and sensitivity operator computation
in the algorithm. Furthermore, in Algorithm \ref{alg:SensInd} the second subscript in sensitivity indices $S_{k,i}$ denotes dependence of the index upon the $i$th data sample. 
\begin{algorithm}[ht!!]
	\caption{Compute the sensitivity indices}
	\begin{algorithmic}[1] \label{alg:SensInd}
        \STATE \verb+% Data sample generation+ 
		\FOR {$i = 1$ to $n_s$}
		\STATE Draw prior sample $\vec{m}_i$
		\STATE Solve the forward equation 
   $\vec{v}(\vec{u}_i,\vec{m}_i,\thth^*_a) = \vec{0}$ for $\vec{u}_i$
		\STATE Synthesize data samples
$\vec{y}_i = \mathcal{O}\vec{u}_i + \vec{\eta}_i(\thth^*_e)$ \hfill \COMMENT{$\mathcal{O}$ observes $\vec{u}$ at measurement locations}
		\ENDFOR
\medskip
        \STATE \verb+% Computation of the Bayes risk sensitivities+
		\FOR {$i = 1$ to $n_s$}
		\STATE Solve the discretized inverse problem for $\vec{m}_i^*(\vec{d}_i,\thth^*)$
		\STATE Solve $-\mat{H}_i\vec{z}_i = \mat{M}(\vec{m}_i^* - \vec{m}_i)$ for $\vec{z}_i$
		\STATE Compute $\vec{r}_i = \mat{B}_i^{\top}\vec{z}_i$
		\ENDFOR 
		\STATE $\DR$ = $\frac{2}{n_s}$$\sum_{i = 1}^{n_s} \vec{r}_i$
		\STATE Compute $\S_k$ and $\S_k^j$ for all $k = 1,\dots,K$ and $j = 1,\dots,n_{\theta}$, see \eqref{eq:Risk_sens}
\medskip
		\STATE \verb+% Computation of the average MAP point sensitivities+
		\FOR {$i = 1$ to $n_s$}
		\STATE Compute $S_{k,i}$ for $k = 1,\dots,K$
			\FOR {$j = 1$ to $n_{\theta}$}
			\STATE Compute $S_{k,i}^j$ for $k = 1,\dots,K$
			\ENDFOR
		\ENDFOR 
		\STATE Compute averaged generalized sensitivities $S_{k} \approx \bar{S}_{k} = \frac{1}{n_s}\sum_{i=1}^{n_s} S_{k,i}$, see \eqref{eq:MAP_sens}
		\STATE Compute averaged pointwise sensitivities $S_k^j \approx  \bar{S}_k^j = \frac{1}{n_s}\sum_{i=1}^{n_s} S_{k,i}^j$, see \eqref{eq:MAP_sens}
	\end{algorithmic}
\end{algorithm}
\FloatBarrier

\subsection{Computational Costs}
Here we discuss the areas of high computational cost in Algorithm~\ref{alg:SensInd}.  
To gain computational efficiency, we rely on some key tools from
PDE-constrained optimization: inexact Newton-CG for MAP estimation,
adjoint methods gradient and Hessian computation, and low-rank approximations for
efficient computation of inverse Hessian 
applies~\cite{Bui-ThanhGhattasMartinEtAl13,PetraMartinStadler2014,
FlathWilcoxHill2011,MartinWilcoxBurstedde2012}. In particular, by combining
methods that make maximum use of the problem structure, we ensure that the
computational complexity of our approach, in the terms of the number of 
PDEs solves, does not scale with the dimension of the discretized inversion parameter.

\textbf{Generate data samples.} 
We solve the forward problem $n_s$ times and use the resulting solutions to generate data. 

\textbf{MAP point solves.} 
We solve the inverse problem $n_s$ times (line 9) using an inexact Newton conjugate gradient line search algorithm with Armijo backtracking. Each Newton step requires 2 PDE solves to compute the gradient and an additional $2I$ PDE solves to compute the Hessian apply where $I$ is the number of iterations required by the CG 
solver to find an appropriate search direction. Thus the total cost is $2L + 2LI$ PDE solves where $L$ is the number of Newton steps taken. This cost in PDE solves multiplied by the number of samples $n_s$ becomes quite significant. However, since the samples drawn from the prior are independent of each other, these computations can be performed in parallel. We also note that we initialize the MAP point solves with the prior samples used to generate data samples.

\textbf{Evaluating inverse Hessian applies.} 
We now address the problem of repeated application of the inverse Hessian, which is required to compute both Bayes risk and MAP point sensitivities in lines 10, 17, and 19. We note that if one only wishes to compute Bayes risk sensitivities, this will not require repeated use of the same Hessian inverse, and line 10 can be evaluated with PCG. Assuming that MAP point sensitivities are also desired, we can offset this cost by computing a low-rank approximation with the Lancoz method to apply the Hessian inverse efficiently, as detailed in \cite{Bui-ThanhGhattasMartinEtAl13}. After computing this low-rank approximation, application of the Hessian inverse can be approximated by matrix-vector products. The computational cost of the Lanczos method is $2r + 2$ PDE solves, where $r$ is the rank of the desired approximation. 

\textbf{Computing Bayes risk sensitivities.}
The sensitivity operator of Bayes risk is a vector, so we built this operator directly before computing indices. We begin this discussion by noting that we can solve the state and adjoint equations around the MAP point once for each data sample, and reuse these solves for each Hessian $\mat{H}_i$ and mixed derivative operator $\mat{B}_i$ or $\mat{B}^{\top}_i$ apply. Each Hessian apply requires 2 additional PDE solves (in addition to the forward and adjoint solves) for the incremental state and incremental adjoint equations. These incremental equation solves can be reused to compute the application of $\mat{B}^{\top}_i$, while $\mat{B}_i$ apply requires 2 more PDE solves for the modified incremental equations. 

\textbf{Computing MAP point sensitivities.} 
The greatest computational cost in estimating the MAP point sensitivities comes in the repeated application of $\DM$ to standard basis vectors $\vec{e}_i$ (line 19) to compute $n_{\theta}$ pointwise sensitivity indices for all $n_s$ data samples. As mentioned previously, this cost is significantly reduced by pre-computing a low-rank approximation that allows for fast Hessian inverse applications. It is also important to note that we reuse the inverse problem solves from computing the Bayes risk sensitivities in computing the MAP point sensitivities and we do not require any additional inverse problem solves here. Due to these various computational savings, we can estimate the MAP point sensitivities through sample averaging at a significantly reduced cost. 

\begin{center}
	\begin{table}[!!!ht]
		\centering
		\caption{Computational costs summary}
		\begin{tabular}
			{ l | l }
			\textbf{Computation} & \textbf{Significant Cost per Sample ($n_s$)}  \\ \hline \hline
			Data Generation & 1 PDE solve \\ \hline
			Inverse Problem Solves & $2L + 2LI$ PDE solves for $L$ Newton steps \\
			& and $I$ PCG iterations \\ \hline
			Hessian Inverse Approximation & $2r + 2$ PDE solves where $r$ is the rank of  \\
			&  the desired approximation \\ \hline
			Bayes risk sensitivities & 2 PDE solves \\ \hline
			MAP point sensitivities & $2n_{\theta}$ PDE solves 
		\end{tabular}
		\label{tb:Costs}
	\end{table}
\end{center}

The discussed computational costs are summarized in Table \ref{tb:Costs} for
clarity.  We remark that for the problem considered in the present work
$n_{\theta}$ is not very large. For problems with a large number of
complementary parameters, computing a suitable low-rank approximation of
$\mat{B}$ may be helpful to reduce the cost of computing many MAP point
sensitivities. We plan to investigate this in our future work.  

\section{Model Problem} \label{sec:ModelProblem}
In this section, we present a model inverse problem, involving heat flow across
a conductive surface, that  will be used to study our HDSA framework.  We begin
by describing the forward problem in Section~\ref{sec:forward} followed by the 
setup of the Bayesian inverse problem in Section~\ref{sec:bayes}.

\subsection{Forward Model}\label{sec:forward}
Consider the problem of infering the log-conductivity field of a
medium from measurements of temperature. Focusing on a cross section, 
we consider the problem in two space dimensions. The forward
problem is governed by the following elliptic PDE, modeling steady state heat
conduction on a unit square domain $\Omega$ with boundary $\partial \Omega =
\cup_{i=1}^{4}\Gamma_i$, where $\Gamma_1$, $\Gamma_2$, $\Gamma_3$, and
$\Gamma_4$ denote the bottom, right, top, and left edges of $\Omega$
respectively,
\begin{subequations} \label{eq:PDE}
	\begin{alignat}{2}
	-\nabla\cdot(e^m\nabla u) &= f && \qquad \text{in } \Omega, \\
	e^m \nabla u \cdot n &= 0 &&\qquad\text{on } \Gamma_1 \cup \Gamma_3,\\
	e^m \nabla u \cdot n &= \beta (T_\text{amb} - u) &&\qquad \text{on } \Gamma_2,\\
	e^m \nabla u \cdot n &= s &&\qquad \text{on } \Gamma_4.
	\end{alignat}
\end{subequations}
In this model, the inversion parameter $m(\vec{x})$ is a function representing the 
log of the heat conductivity of the non-homogeneous two-dimensional surface. 
We let $u(\vec{x})$ denote the temperature, $f(\vec{x})$ the heat source in the domain, $\beta$ the heat transfer coefficient of the medium, $T_{\text{amb}}$ the ambient temperature of the medium, and $s(x_2)$ a boundary heat source function representing heat entering the domain from the left boundary. In this model problem, the equations in \eqref{eq:PDE} are dimensionless and we let $T_{\text{amb}} = 22$ and consider the heat transfer coefficient $\beta$ to be an uncertain auxiliary parameter with a nominal value of $\beta = 1$. 

\begin{figure} [ht!]
	\centering
	\includegraphics[width=0.48\linewidth]{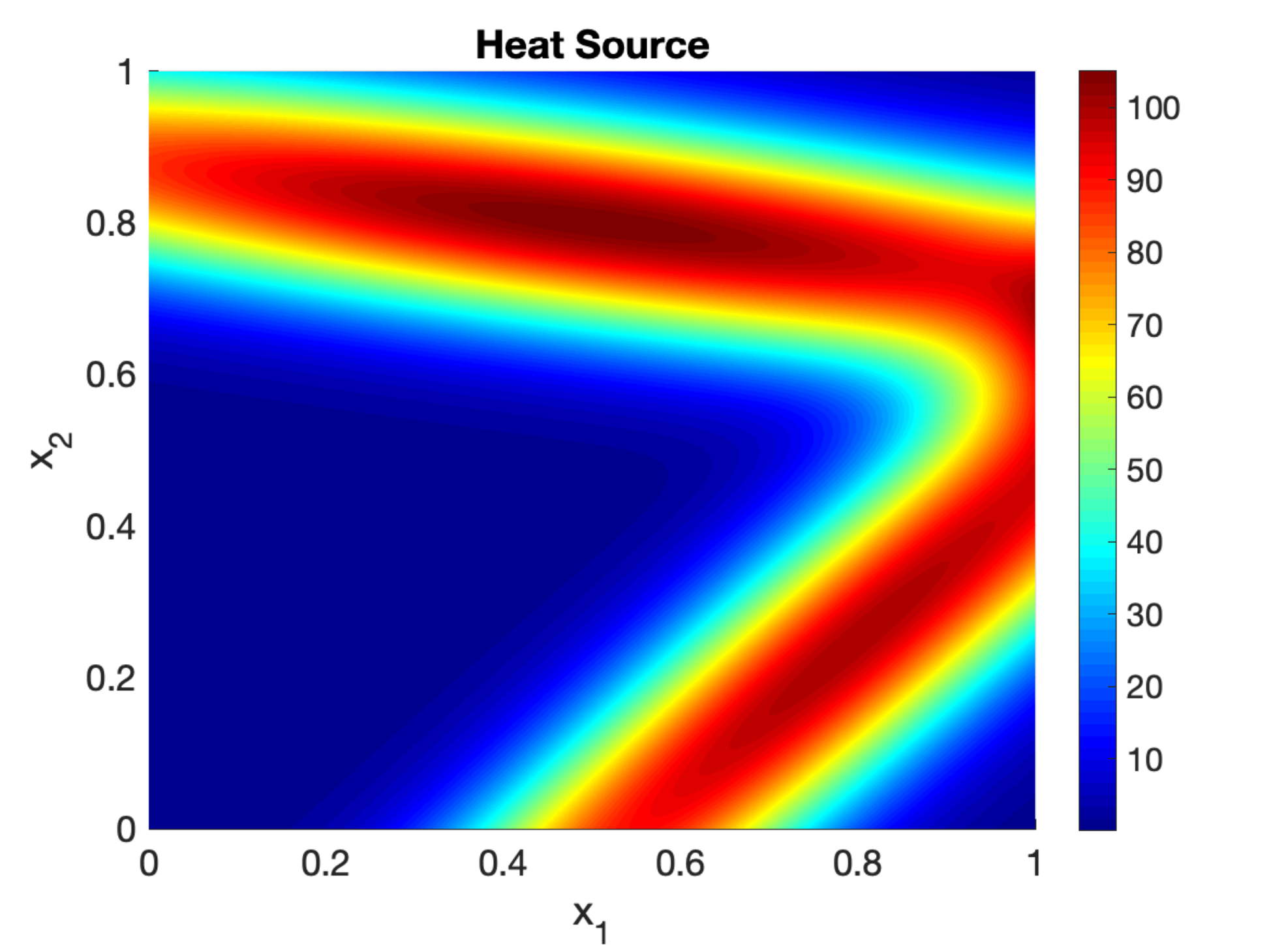}
	\caption{The heat source function $f(\vec{x})$}
	\label{fig:truesource}
\end{figure}

The boundary heat source $s(x_2)$ is modeled as follows,
\[
s(x_2) = s_1\exp\left(-\left(\frac{x_2-s_3}{s_2}\right)^2\right)
\]
with auxiliary parameters $s_1, s_2, \text{and } s_3$ fixed at nominal values $s_1 = 30, s_2 = .1, \text{and } s_3 = .65$. The auxiliary parameters consist of the amplitude, spread, and location of the boundary heat source respectively. The heat source in the domain $f(\vec{x})$ is modeled as,
\[
f(\vec{x}) = f_1\exp\left[-\frac12(\vec{x} - \vec{w})^{\top} \mat{C}_1 (\vec{x} - \vec{w}) \right] + 
f_2\exp\left[-\frac12(\vec{x} - \vec{z})^{\top} \mat{C}_2 (\vec{x} - \vec{z}) \right], \quad \text{with}
\]
\[
\mat{C}_1 = \begin{bmatrix}
\frac{\cos^2(\gamma_1)}{\sigma_{x_1}^2} + \frac{\sin^2(\gamma_1)}{\sigma_{x_2}^2} & \frac{\sin(2\gamma_1)}{2\sigma_{x_2}^2} - \frac{\sin(2\gamma_1)}{2\sigma_{x_1}^2} \\
\frac{\sin(2\gamma_1)}{2\sigma_{x_2}^2} - \frac{\sin(2\gamma_1)}{2\sigma_{x_1}^2}  &
\frac{\sin^2(\gamma_1)}{\sigma_{x_1}^2} + \frac{\cos^2(\gamma_1)}{\sigma_{x_2}^2}
\end{bmatrix} \text{ and }
\mat{C}_2 = \begin{bmatrix}
\frac{\cos^2(\gamma_2)}{\sigma_{x_1}^2} + \frac{\sin^2(\gamma_2)}{\sigma_{x_2}^2} & \frac{\sin(2\gamma_2)}{2\sigma_{x_2}^2} - \frac{\sin(2\gamma_2)}{2\sigma_{x_1}^2} \\
\frac{\sin(2\gamma_2)}{2\sigma_{x_2}^2} - \frac{\sin(2\gamma_2)}{2\sigma_{x_1}^2} & \frac{\sin^2(\gamma_2)}{\sigma_{x_1}^2} + \frac{\cos^2(\gamma_2)}{\sigma_{x_2}^2} 
\end{bmatrix}.
\]
In this formulation, $f_1$ and $f_2$ control the amplitude of the heat sources, $\vec{w}$ and $\vec{z}$ control the centers of the two sources, $\gamma_1$ and $\gamma_2$ their respective tilt angles, and $\sigma_{x_1}$ and $\sigma_{x_2}$ the spread of the heat sources in the $x_1$ and $x_2$ directions respectively. For this problem we fix these parameters at the following nominal values: $f_1 = 100, f_2 = 105, \vec{w} = (.8,.25), \vec{z} = (.5,.8), \gamma_1 = -\pi/4, \gamma_2 = .15, \sigma_{x_1} = .8, \sigma_{x_2} = .1$. We consider the amplitude, center point, and angle of each bar to be uncertain and thus let $f_1,f_5,w_1,w_2,z_1,z_2,\gamma_1,$ and $\gamma_2$ be the auxiliary parameters for the right hand side heat source $f(\vec{x})$. Figure \ref{fig:truesource} depicts this heat source in the domain. 

\FloatBarrier

We note that this model problem has been kept intentionally simple to
aid in the interpretation and understanding of the complicated
algorithmic methodology. Even so, this example is motivated by many
uncertainties surrounding additive manufacturing processes (such as
powered bed laser fusion) that cause high residual stresses and even
defects in final parts.  Variability in the powder material, boundary
conditions, rasterization patterns, and laser power result in uneven
heat distribution with problematic micro-crystallographic structures
and inhomogeneous material properties.  Although the underlying
physics for additive manufacturing is more complicated, our model
problem conceptually demonstrates the ability of our approach to
provide insight into a complicated application area.

\subsection{Prior Measure and State Solution}
\label{sec:bayes}
In many inverse problems a ``true solution" is chosen to synthesize data and evaluate the accuracy of the proposed methodology. Note that we do not have any such ``true solution" here and instead we compute data from samples of the prior distribution. 
We specify the Bayesian prior on $m$ as a Gaussian random field on $\Omega$ with mean $m_{pr}$ and covariance operator $\Cpr$. We model the prior mean as a sinusoidal function: 
\[
m_{pr}(\vec{x}) = 1.5\sin(2\pi x_1)\cos(2\pi x_2) + 2.
\]
We let the covariance operator $\Cpr$ be the inverse of a squared elliptic differential operator $\mathcal{A}$, where $m = \mathcal{A}^{-1}s$ satisfies
\[
\alpha\int_{\Omega} (\Phi \nabla m) \cdot \nabla q + mq \ d\vec{x} = \int_{\Omega} sq \ d\vec{x}
\]
for all $q \in H^1(\Omega)$, with $\alpha = 5$, and $\Phi = .01$. This formulation of the prior covariance ensures that $\Cpr$ is trace class and provides a computationally convenient formulation. For more details see \cite{Bui-ThanhGhattasMartinEtAl13}. 

\begin{figure}[ht!]
	\centering
	\includegraphics[width=0.45\linewidth]{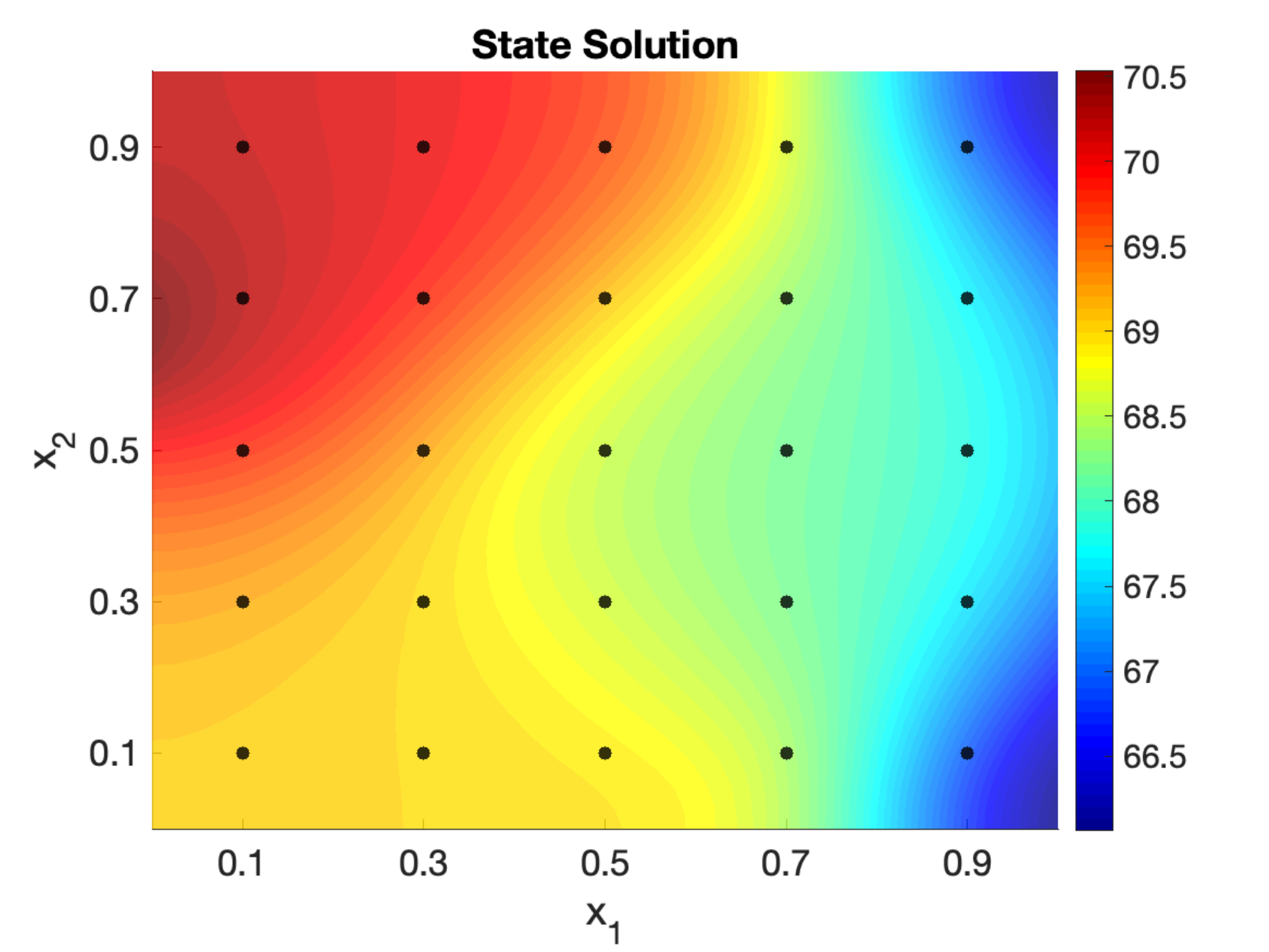}
	\caption{The state solution of the governing system of partial differential equations with the experimental sensor locations indicated by filled black circles.}
	\label{fig:statesolution}
\end{figure}

Measurements are collected on an evenly spaced 5$\times$5 grid of observation locations depicted in Figure \ref{fig:statesolution}. We consider the standard deviation of the noise in each data measurement to be our uncertain experimental parameters. Additive Gaussian noise models ``error" in our data and we assume the measurements are uncorrelated, with nominal standard deviations of $\sigma = .1$, thus $\mat{\Gamman} = \sigma^2\mat{I}$. Although we allow the measurement standard deviations to take the same nominal value, we consider each standard deviation individually when computing sensitivities of the solution. Perturbing the noise standard deviation will also result in a perturbation of the noise realization $\eta_i \sim \mathcal{N}(0, \Gamman)$, directly proportional to the multiplicative perturbation of $\sigma_i$. Therefore the experimental parameters $\thth_e$ enter the inverse problem through the cost function~\eqref{eq:MAP_cost}, both in the noise covariance matrix $\Gamman(\thth_e)$ and the data measurements $y(\thth_e)$ which depend on the noise realizations. 

The solution of the governing PDE system detailed in \eqref{eq:PDE} at the nominal parameter values with $m$ fixed at the prior mean is depicted in Figure \ref{fig:statesolution}. 

\FloatBarrier

\section{Results} \label{sec:Results}
\begin{figure} [ht!]
	\centering
	\includegraphics[trim={200 0 200 0},clip,width=0.9\linewidth]{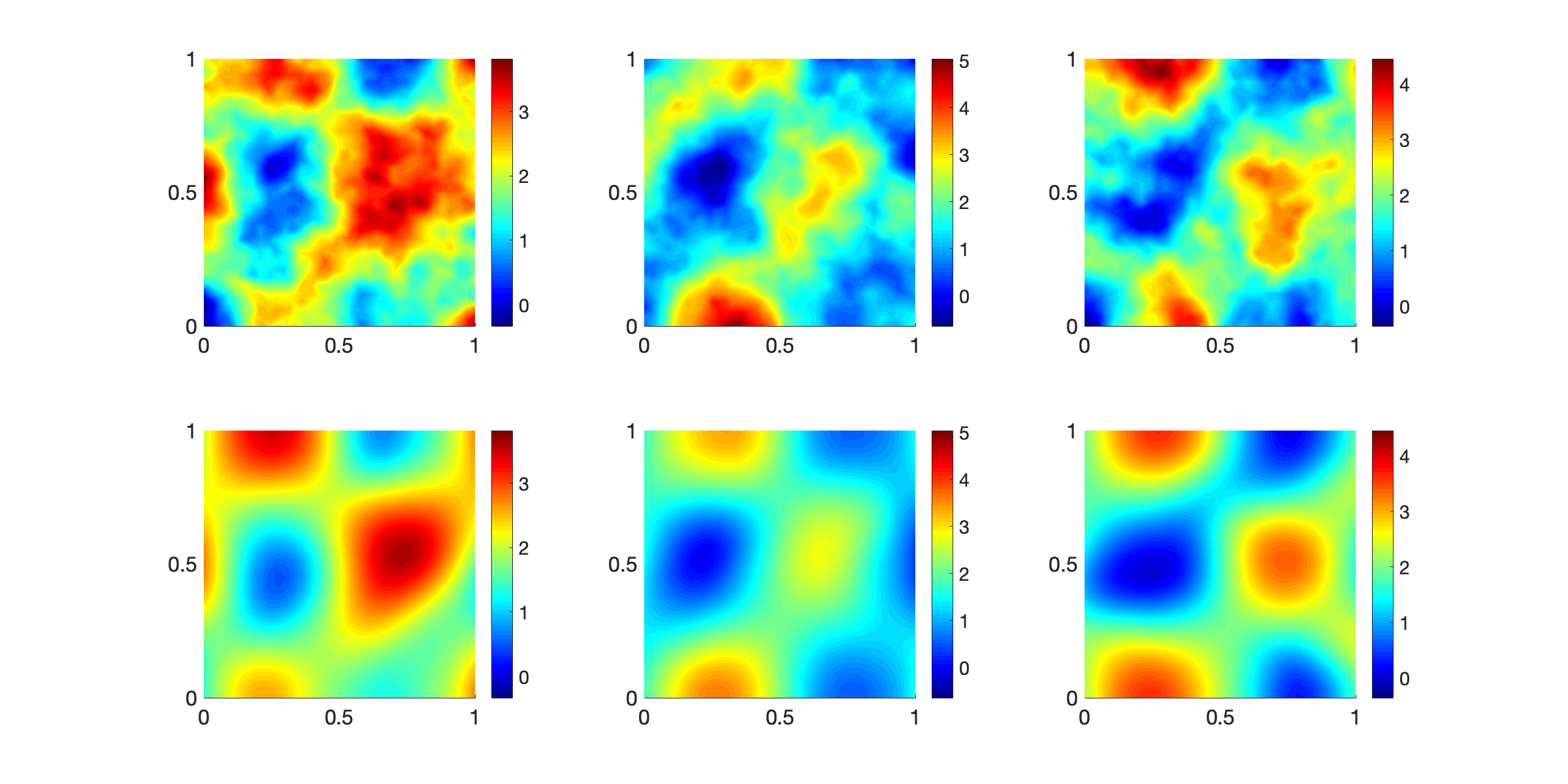}
	\caption{Top: three samples from the log-conductivity prior. Bottom: inverse problem MAP point estimates solved using data generated by the above prior sample.}
	\label{fig:sampmaprecon2}
\end{figure}

Using the model of heat flow across a conductive surface from Section \ref{sec:ModelProblem}, we solve an inverse problem to estimate the posterior distribution. Following Algorithm \ref{alg:SensInd} to evaluate our Bayesian hyper-differential sensitivities, we take samples from the prior distribution on $m$ and push them through the forward mapping to generate noisy data. Each data sample is then used to solve \eqref{eq:MAP_opt}, giving a unique MAP point reconstruction for each sample. To illustrate this process, we present three prior samples and their corresponding MAP point reconstructions in Figure \ref{fig:sampmaprecon2}. 
\FloatBarrier

Each MAP point is attempting to estimate the above prior sample from noisy data. This example is illustrative in that it gives us some insight into Bayes risk, which measures the average difference in norm between the prior samples (top) and the inferred MAP points (bottom). 

In Section \ref{sec:complementary_parameters} we detail how perturbations of the complementary parameters are modeled. Following this we present and discuss the significance of the generalized sensitivities of the complementary parameters as well as the pointwise sensitivities of the experimental parameters with respect to Bayes risk (Section \ref{subsec:RiskResults}) and the MAP point (Section \ref{subsec:MAPResults}). We note that the pointwise sensitivities of the auxiliary parameters are identical to their generalized sensitivities as each auxiliary parameter is scalar valued in this model problem. 

\subsection{Modeling Parameter Perturbations}
\label{sec:complementary_parameters}
Suppose $\rho$ is an uncertain scalar parameter. We model our uncertainty in $\rho$ as,
\begin{equation} \label{eq:perturb1}
\rho = \tilde{\rho}(1 + a\theta)
\end{equation}
where $\tilde{\rho}$ is the nominal value, $a$ is a scaling coefficient quantifying our degree of uncertainty, and $\theta \in [-1,1]$ defines a perturbation of $\tilde{\rho}$. Perturbations of vector valued complementary parameters, such as data measurements, are modeled as componentwise scalar perturbations as in \eqref{eq:perturb1}. 

In this particular model problem we use a perturbation scaling coefficient of $a = .05$ for each auxiliary parameter, which represents our uncertainty in that parameter's estimate being $5\%$ of the parameter's nominal value. For the experimental parameters we instead use a scaling coefficient of $a = 1$ to represent that our uncertainty in the standard deviation of the data noise is the full quantity of the standard deviation. 

\subsection{Sensitivities of Bayes Risk} \label{subsec:RiskResults}

\begin{figure}[ht!]
	\centering
	\includegraphics[width=0.9\linewidth]{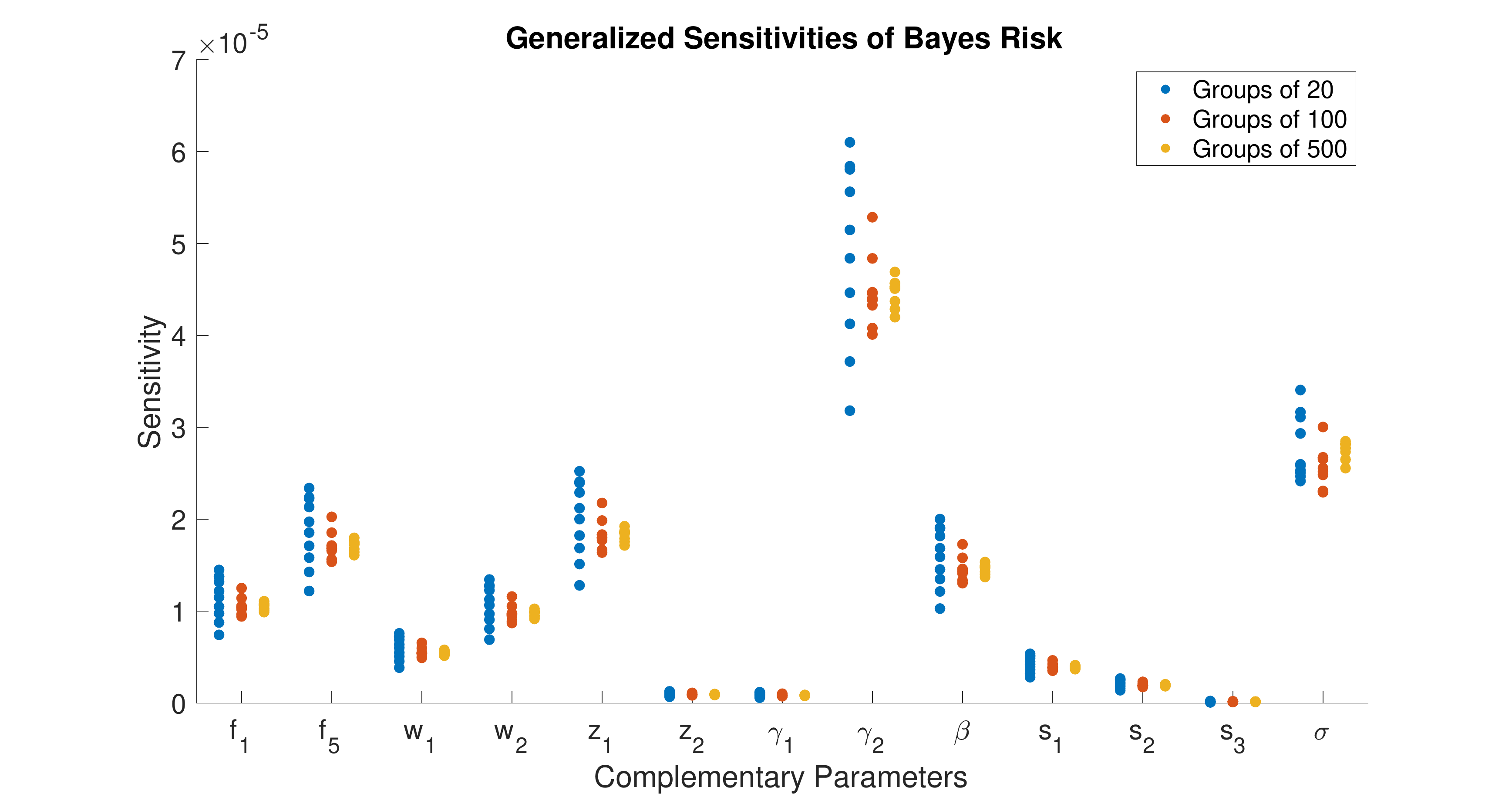}
	\caption{Generalized sensitivities of Bayes Risk to complementary parameters computed for 10 groups of samples of size 20, 100, and 500.}
	\label{fig:GenSensDR}
\end{figure}

The approximate Bayes risk is computed as a sample average as detailed in \eqref{eq:BayesRiskDis}. We present the generalized sensitivities of each complementary parameter with respect to Bayes risk in Figure \ref{fig:GenSensDR}. We study the effect of the sample size on the computed sensitivities by comparing generalized sensitivities for Bayes risk computed from ten groups of 20 samples, ten groups of 100 samples, and ten groups of 500 samples, each taken randomly from a group of 3000 pre-computed samples. 

\FloatBarrier

First let us discuss the spread of the samples. 
We can see that while the groups of 100 samples can sometimes vary significantly, such as in the case of $\gamma_2$, they generally capture the rankings of the parameters relative to one another correctly. Thus, if we are primarily concerned with determining the relative importance of the parameters compared to each other we may conclude that 100 samples provides sensitivity estimates that suit our needs. 

Next, we note that Bayes risk is most sensitive to the tilt angle of the second domain heat source $\gamma_2$. We observe that as the tilt angles $\gamma_1$ and $\gamma_2$ are changed, they can overlap in the domain interior, causing a large increase in heat where the overlap occurs and will result a significant change in $f$. One possible reason $\gamma_2$ is so important is that even a relatively small perturbation will result in increased or decreased overlap of these bars in the domain. Of secondary importance are the heat amplitude ($f_5$) and center in the $x_1$ direction ($z_1$) of the second domain heat source, heat transfer coefficient ($\beta$), and the standard deviations of data noise ($\vec{\sigma}$). This sensitivity information can then be used by an experimenter to inform their experimental design choices for this problem. 
To accurately estimate the Bayes risk for this problem as a measure of posterior uncertainty, it is more important to invest resources in ensuring that the parameters $\gamma_2, f_5, z_1, \beta$, and $\vec{\sigma}$ are accurately estimated than the other complementary parameters. Specifically, we can interpret these sensitivities as ``a 5\% perturbation in the scalar auxiliary parameters or a norm-1 perturbation in the experimental parameters ($\vec{\sigma}$) will result in a perturbation of Bayes risk proportional to the sensitivity." 

While these sensitivities appear to be very small, we note that the problem is highly diffusive and steady state. Both of these factors are likely making the problem highly insensitive to perturbations of complementary parameters. This in of itself showcases the benefits of using HDSA. For such an insensitive problem, it would be extremely difficult to gather any kind of intuition or conclusion as to the relative importance of various parameters a priori. With our framework however, we can rigorously determine the relative importance of uncertain parameters before any physical experimentation is done, even for highly insensitive problems, which is valuable to experimenters who seek to efficiently allocate experimental resources. 

\begin{figure}[ht!]
	\centering
	\includegraphics[width=0.65\linewidth]{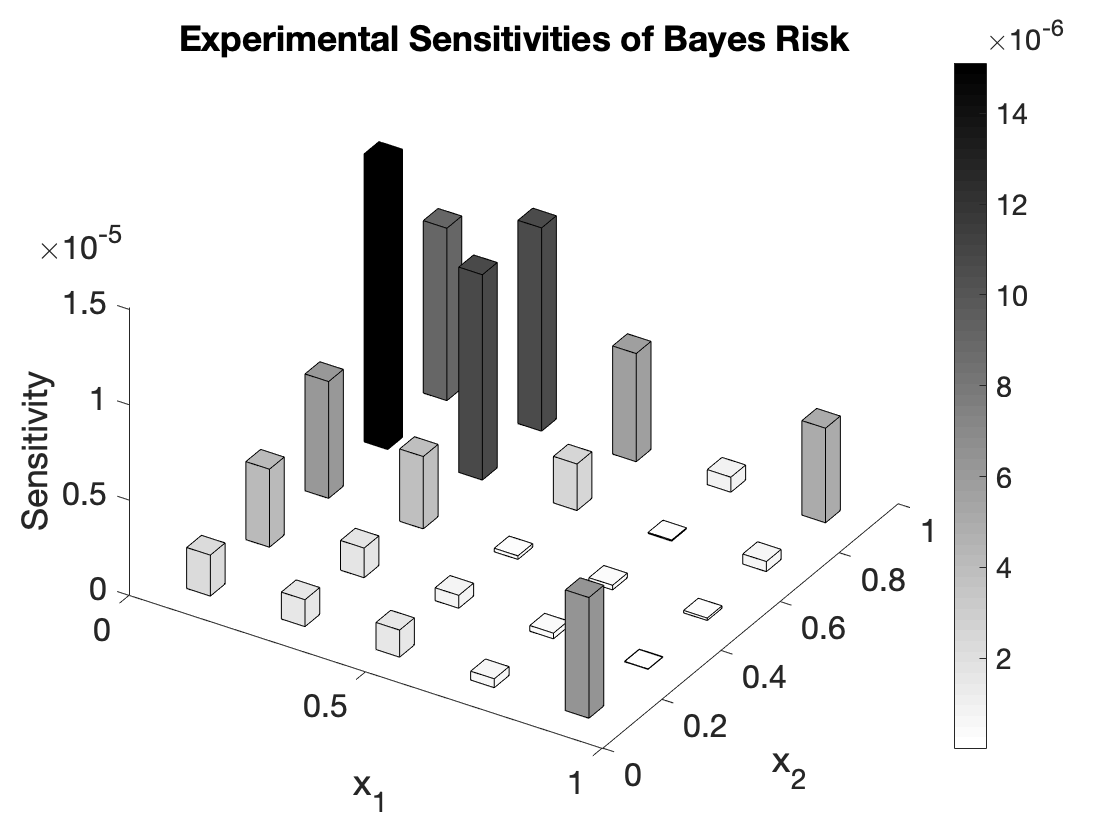}
	\caption{Point-wise sensitivities of Bayes Risk to experimental parameters computed from 1000 samples. Each bar represents the sensitivity of perturbing the standard deviation of the data noise at that particular sensor location.}
	\label{fig:ExpSenseDR}
\end{figure}

Next we study the pointwise sensitivities of Bayes risk to the experimental parameters,  the standard deviation of noise in the data measurements, presented in Figure \ref{fig:ExpSenseDR}. By perturbing the noise, we model perturbations of each collected data measurement in a way that we can experimentally control through sensor accuracy.

\FloatBarrier

We first note the scale of sensitivities presented in Figure \ref{fig:ExpSenseDR}. Although these sensitivities are very small (on the order of $10^{-5}$ or $10^{-6}$), this is not entirely unexpected given the scale of the generalized sensitivities presented in Figure \ref{fig:GenSensDR}. Indeed, we would expect that perturbing a single data measurement's noise would not result in a very large change in Bayes risk. We can see that the sensors grouped around small values of $x_1$ and large values of $x_2$ are most important with respect to Bayes risk. Thus, we can conclude that the data measured at these sensors is the most important to collect accurately for the purposes of estimating our measure of posterior uncertainty. We observe that these sensors are located in the region that the state solution depicted in Figure \ref{fig:statesolution} is largest. This is also the region near the boundary source term $s$. We also note that the sensors at (.9,.1) and (.9,.9) are relatively important, which are located in the areas where the state solution is smallest. These results provide information that may not be obvious a priori and helps practitioners understand what parameters and sensor measurements the solution is most sensitive to. 

\subsection{Sensitivities of the MAP Point} \label{subsec:MAPResults}
\begin{figure}[ht!]
	\centering
	\includegraphics[width=0.9\linewidth]{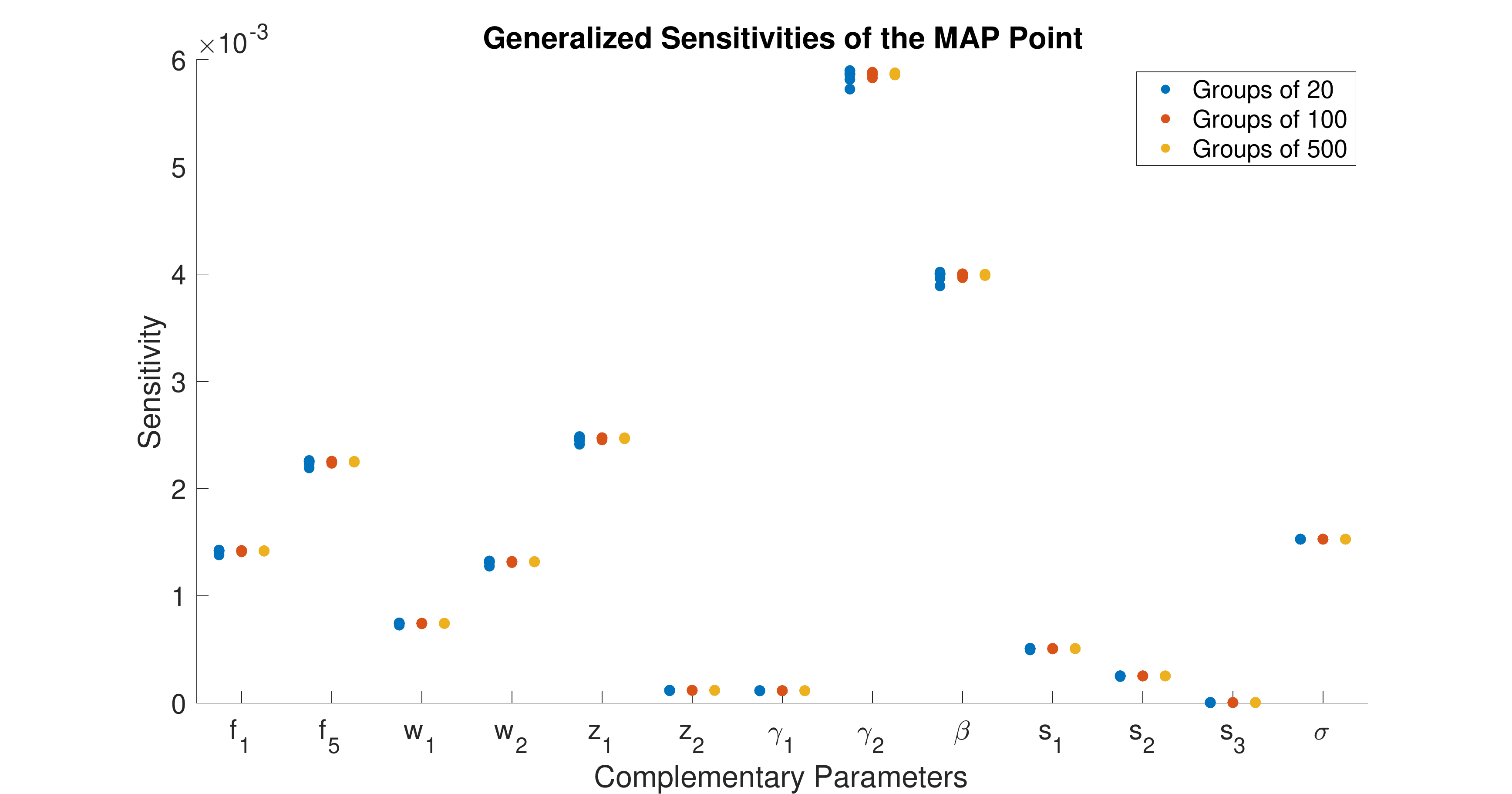}
	\caption{Averaged generalized sensitivities of the MAP point to complementary parameters computed for 10 groups of sensitivities of size 20, 100, and 500.}
	\label{fig:GenSensD}
\end{figure}

We now study the averaged generalized sensitivities of the MAP point. As was done previously, we study the effect of the sample size on the generalized sensitivities. This is done by computing generalized MAP point sensitivities for 3000 data samples. We then randomly select and average ten groups of 20 sensitivities, ten groups of 100 sensitivities, and ten groups of 500 sensitivities, which are plotted in Figure \ref{fig:GenSensD}. 

\FloatBarrier

In this case we can see that even groups of 20 sensitivities produce little variation in the averaged sensitivity measure. Thus, we can conclude that for this application, using a sample average of just 20 sensitivities provides sufficient accuracy for our purposes. Furthermore, we notice that the generalized sensitivities of the MAP point are significantly greater in magnitude than those computed for Bayes risk. For this problem, it appears that the MAP point is more sensitive to perturbations in the complementary parameters than the posterior uncertainty is. We see that the MAP point has greatest sensitivity to $\gamma_2, \beta, z_1$, and $f_5$. It is interesting to note that for Bayes risk, $\vec{\sigma}$ and $\beta$ had the second and fifth greatest sensitivity, respectively. In contrast, the sensitivity rankings of these two parameters have switched places with respect to the MAP point. 

\begin{figure}[ht!]
	\centering
	\includegraphics[width=0.7\linewidth]{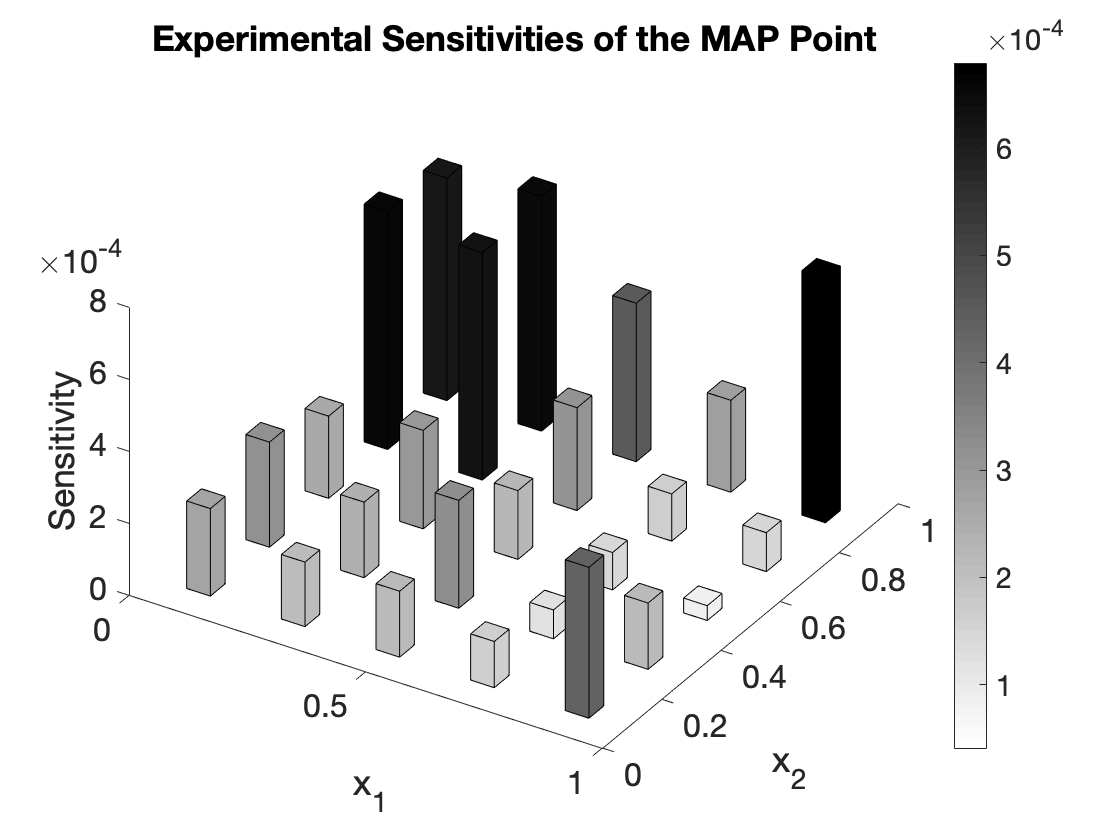}
	\caption{Point-wise sensitivities of the MAP point to experimental parameters computed as an average of 20 sample sensitivities.}
	\label{fig:ExpSenseD}
\end{figure}

Finally, we examine the pointwise sensitivities of the MAP point to the experimental parameters depicted in Figure \ref{fig:ExpSenseD}. Each pointwise sensitivity is computed as an average of 20 sensitivities computed from different data samples. We compared these pointwise sensitivities with those computed from an average of 1000 sensitivities, and as our study on sample size in Figure \ref{fig:GenSensD} would indicate, there was minimal difference. 

\FloatBarrier

We can again see that the sensors grouped around small values of $x_1$ and large values of $x_2$ are most important with respect to the MAP point. We observe that for this problem the sensors with greatest importance to the MAP point coincide closely with those sensors that are important for Bayes risk. 

\section{Conclusion} \label{sec:Conclusion}

In this article we take foundational steps in applying hyper-differential
sensitivity analysis (HDSA) to large-scale nonlinear Bayesian inverse problems.
In particular, we focus on HDSA of the MAP point and Bayes risk to the
auxiliary and experimental parameters and present efficient methods for
computing the corresponding HDSA indices.  Performing HDSA is important as it reveals the auxiliary
parameters the inverse problem is most sensitive to. Moreover, HDSA with respect
to measurement data helps identify the measurements that are important to the
solution of the inverse problem, and can guide design of experiments by investing
resources to obtain good quality data from important measurement points.

It is also important to note that the Bayesian formulation allows for the
computation of HDSA indices prior to conducting
experiments. Namely, we use the information encoded in the Bayesian inverse
problem to obtain likely realizations of measurement data, 
which are used to compute the Bayes risk
sensitivities and average MAP point sensitivities.  This is a key factor
that makes this approach attractive for HDSA of Bayesian inverse
problems, while minimizing experimental costs. 

While the steady state heat conduction model presented in Section
\ref{sec:ModelProblem} is an academic model problem, it has many features that
are seen in real applications.  We found that the tilt angle, heat amplitude,
and center in the horizontal direction of volume heat source as well as the
heat transfer coefficient and data noise were the parameters that both the
Bayes risk and the MAP point were most sensitive to. We also determined which
sensors provide the most informative data and found that for this problem the
Bayes risk is generally less sensitive to perturbations of the complementary
parameters than the MAP point is.  Such observations can be instrumental in
areas such as additive manufacturing.  By applying the proposed methods to
additive manufacturing problems, one can determine a priori which experimental
factors the inverse problem solution will be most sensitive to and thereby
guide the calibration of equipment tolerances with this information.

The MAP point is a key point estimator for the inversion parameters and
performing HDSA on this quantity provides valuable insight regarding the
sensitivity of the inverse problem to complementary parameters. On the other
hand, Bayes risk provides a measure of the statistical quality of the estimated
parameters, and is a common utility function in decision theory.  Additionally,
up to a linearization, Bayes risk can be considered as a proxy for posterior
uncertainty. These considerations, coupled with the fact that the methods for
HDSA of Bayes risk build on methods for HDSA of MAP point, made Bayes risk a
suitable HDSA QoI in first steps towards HDSA of Bayesian inverse problems.

In our future work, we plan to investigate HDSA of different quantities such as
average posterior variance or expected information gain. Suitable
approximations of the posterior, such as a Laplace approximation, can be
considered, to mitigate the high cost of HDSA of such quantities in large-scale
nonlinear inverse problems. Another interesting line of inquiry is to use HDSA
within the context of optimal experimental design
(OED) under uncertainty~\cite{AlexanderianPetraStadlerSunseri21,
KovalAlexanderianStadler19}. HDSA can reveal model uncertainties that the OED
criterion is most sensitive to and thus must be accounted for in the optimal
design process. On the other hand, model uncertainties the design criterion is
less sensitive to may be fixed at some nominal values, hence reducing the
complexity of OED under uncertainty problem. 

\section*{Acknowledgements}
This paper describes objective technical results and analysis. Any subjective views or opinions that might be expressed in the paper do not necessarily represent the views of the U.S. Department of Energy or the United States Government. Sandia National Laboratories is a multimission laboratory managed and operated by National Technology and Engineering Solutions of Sandia LLC, a wholly owned subsidiary of Honeywell International, Inc., for the U.S. Department of Energy's National Nuclear Security Administration under contract DE-NA-0003525. SAND2022-1279 O. This work was supported by the US Department of Energy, Office of Advanced Scientific Computing Research, Field Work Proposal 20-023231.

The work of I. Sunseri and A. Alexanderian was supported in part by the
National Science Foundation under grant DMS-1745654. 
Additionally, the work of A. Alexanderian was also supported in part by the
National Science Foundation under grant DMS-2111044. 

\section*{Bibliography}
\bibliographystyle{unsrt}
\bibliography{refs}

\end{document}